\title
{A short proof of Kneser's addition theorem \\ for abelian groups}
\author{
	Matt DeVos	\\
	{\normalsize mdevos@sfu.ca}
}

\date{}

\documentclass[12pt]{article}

\usepackage{amsmath}
\usepackage{amsfonts}

\begin{document}
\bibliographystyle{plain}
\maketitle
\setcounter{page}{1}
\newtheorem{theorem}{Theorem}
\newtheorem{lemma}[theorem]{Lemma}
\newtheorem{corollary}[theorem]{Corollary}
\newtheorem{proposition}[theorem]{Proposition}
\newtheorem{definition}[theorem]{Definition}
\newtheorem{claim}{Claim}
\newtheorem{conjecture}[theorem]{Conjecture}
\newtheorem{observation}[theorem]{Observation}
\newtheorem{problem}[theorem]{Problem}
\newtheorem{question}[theorem]{Question}

\begin{abstract}
Martin Kneser proved the following addition theorem for every abelian
group $G$.  If $A,B \subseteq G$ are finite and nonempty, then 
$|A+B| \ge |A+K| + |B+K| - |K|$ where 
$K = \{ g \in G \mid g+A+B = A+B \}$.  Here we give a short proof of this 
based on a simple intersection union argument.
\end{abstract}

Throughout we shall assume that $G$ is an additive abelian 
group.  If $A,B \subseteq G$ and $g \in G$, then 
$A + B = \{ a+b \mid \mbox{$a \in A$ and $b \in B$} \}$ and
$A + g = g + A = \{ a+g \mid a \in A \}$.  We define the {\it stabilizer}
of $A$ to be ${\mathcal S}(A) = \{ g \in G \mid A + g = A \}$.  Note
that ${\mathcal S}(A) \le G$.  The goal of this paper is to provide a
short proof of the following theorem.

\begin{theorem}[Kneser \cite{kneser}]
If $A,B \subseteq G$ are finite and nonempty and $K = {\mathcal S}(A+B)$, 
then 
\[ |A + B| \ge |A + K| + |B + K| - |K|.\]
\end{theorem}

\noindent{\bf Proof.}
We proceed by induction on $|A+B| + |A|$.  Suppose that $K \neq \{0\}$ 
and let $\phi : G \rightarrow G / K$ be the canonical homomorphism.  
Then ${\mathcal S}( \phi(A+B) )$ is trivial, so by applying induction to
$\phi(A),\phi(B)$ we have 
\[ |A+B| = 
	|K|(|\phi(A) + \phi(B)|) \ge |K|(|\phi(A)| + |\phi(B)| - 1) 
	= |A+K| + |B+K| - |K|.\] 
Thus, we may assume $K = \{0\}$.  If $|A| = 1$, then the result
is trivial, so we may assume $|A| > 1$ and choose distinct $a,a' \in A$.  
Since $a' - a \not\in {\mathcal S}(B) \subseteq {\mathcal S}(A+B) = \{0\}$, 
we may choose $b \in B$ so that $b+a'-a \not\in B$.  Now by replacing 
$B$ by $B - b + a$ we may assume  $\emptyset \neq A \cap B \neq A$.

Let $C \subseteq A+B$ and let $H = {\mathcal S}(C)$.  We call $C$ 
a {\it convergent} if
\[ |C| + |H| \ge |A \cap B| + |(A \cup B) + H|. \]
Set $C_0 = (A \cap B) + (A \cup B)$ and observe that $C_0 \subseteq A + B$.  
Since $0 < |A \cap B| < |A|$, we may apply induction to $A \cap B$ and 
$A \cup B$ to conclude that $C_0$ is a convergent.  Thus a convergent exists, and 
we may now choose a convergent $C$ with $H = {\mathcal S}(C)$ minimal.
If $H = \{0\}$ then 
$|A+B| \ge |C| \ge |A \cap B| + |A \cup B| - |\{0\}| = |A| + |B| - 1$ 
and we are finished.  So, we may assume $H \neq \{0\}$ (and proceed toward a contradiction).
Since ${\mathcal S}(A+B) = \{0\}$ and ${\mathcal S}(C)=H$, 
we may choose $a \in A$ and $b \in B$ so that 
$a + b + H \not\subseteq A+B$.  Let $A_1 = A \cap (a + H)$, 
$A_2 = A \cap (b + H)$, $B_1 = B \cap (b + H)$, and 
$B_2 = B \cap (a + H)$ and note that $A_1,B_1 \neq \emptyset$.  For $i=1,2$
let $C_i = C \cup (A_i + B_i)$ and let $H_i = {\mathcal S}(A_i + B_i)$.  
Observe that if $A_i,B_i \neq \emptyset$, then $H_i = {\mathcal S}(C_i) < H$.
The following equation holds for $i=1$, and it also holds for $i=2$ if 
$A_2,B_2 \neq \emptyset$.  It follows from the fact that $C_i$ is not a 
convergent (by the minimality of $H$), and induction applied to $A_i,B_i$.  
\begin{eqnarray}
\label{eq1}
|(A \cup B)+H| - |(A \cup B)+H_i| 
 & < & (|C| + |H| - |A \cap B|) - (|C_i| + |H_i| - |A \cap B|)	\nonumber\\
 & = & |H| - |A_i + B_i| - |H_i|	\nonumber\\
 &\le& |H| - |A_i + H_i| - |B_i + H_i|
\end{eqnarray}
If $B_2 = \emptyset$, then 
$|(A \cup B)+H| - |(A \cup B) + H_1| \ge |(a+H) \setminus (A_1 + H_1)| = |H| - |A_1 + H_1|$  
contradicts equation \ref{eq1} for $i=1$.  We get a similar contradiction 
under the assumption that $A_2 = \emptyset$.  Thus $A_2,B_2 \neq \emptyset$ 
and equation \ref{eq1} holds for $i=1,2$.  
If $a + H = b + H$, then $A_1=A_2$ and $B_1=B_2$ and we have
$|(A \cup B) + H| - |(A \cup B) + H_1| \ge
|(a+H) \setminus ((A_1 \cup B_1) + H_1)| \ge 
|H| - |A_1 + H_1| - |B_1 + H_1|$ which contradicts equation \ref{eq1}.  
Therefore, $a + H \neq b + H$.  Our next inequality follows
from the observation that the left hand side of equation \ref{eq1} is
nonnegative, and all terms on the right hand side are multiples of $|H_i|$.
\begin{equation}
\label{eq2}
|H| \ge |A_i| + |B_i| + |H_i|
\end{equation}
Let $S = (a + H) \setminus (A_1 \cup B_2)$ and  
$T = (b + H) \setminus (A_2 \cup B_1)$, and note that $S$ and $T$ are
disjoint.  The next equation follows from the fact that 
$A+B$ is not a convergent (by the minimality of $H$), and
induction applied to $A_i,B_i$.
\begin{eqnarray}
\label{eq3}
 |H| 
 &\ge& |(A \cup B)+H| + |A \cap B| - |C|	\nonumber\\
 &\ge& |S| + |T| + |A \cup B| + |A \cap B| - |A+B| + |A_i + B_i| \nonumber\\
 & > & |S| + |T| + |A_i| + |B_i| - |H_i|	
\end{eqnarray}
Summing the four inequalities obtained by taking equations \ref{eq2} and 
\ref{eq3} for $i=1,2$ and then dividing by two yields
$2|H| > |A_1| + |B_2| + |S| + |A_2| + |B_1| + |T|$.  However, 
$a + H = S \cup A_1 \cup B_2$
and $b + H = T \cup A_2 \cup B_1$.  This final contradiction 
completes the proof.
\quad\quad$\Box$

\end{document}